\documentclass[a4paper, 12pt]{article}

\usepackage{kotex}
\usepackage{amsfonts}
\usepackage {amssymb}
\usepackage {amsmath}
\usepackage {amsthm}
\usepackage{graphicx}
\usepackage{multirow}
\usepackage{xypic}
\usepackage {amscd}
\usepackage{mathrsfs}
\usepackage[colorlinks, linkcolor=blue, anchorcolor=black, citecolor=red]{hyperref}
\usepackage{enumerate}
\usepackage{enumitem}
\usepackage{geometry}  
\usepackage{tikz-cd}
\usepackage[all]{xy}

\newcommand{\oX}{X}
\newcommand{\bC}{\mathbb{C}}

\newcommand{\bZ}{\mathbb{Z}}
\newcommand{\bR}{\mathbb{R}}
\newcommand{\bCP}{\mathbb{CP}^n}
\newcommand{\tCP}{\mathbb{CP}^n\times\mathbb{CP}^n}
\newcommand{\bP}{\mathbb{P}}

\newcommand{\cO}{\mathcal{O}}

\newcommand{\cN}{\mathcal{N}}

\newtheorem{thm}{Theorem}[section]
\newtheorem{prop}[thm]{Proposition}
\newtheorem{defn}[thm]{Definition}

\newtheorem{rem}[thm]{Remark}

\newtheorem{lem}[thm]{Lemma}

\newtheorem{defn-prop}[thm]{Definition-Proposition}

\begin{document}

\title{Zoll Manifolds of Type $\bCP$ with Entire Grauert Tubes}
\author{Chi Li, Kyobeom Song}
\date{}

\maketitle

\abstract{
We show that a Zoll manifold of type $\bCP$ with an entire Grauert tube is isometric to $\bCP$ with the canonical Fubini--Study metric, up to constant multiplication.
  }

\tableofcontents

\section{Introduction and main results}

A Riemannian manifold $(M,g)$ is called a \emph{Zoll manifold} if there exists a positive real number $l$ satisfying:
\begin{itemize}
\item Every geodesic $\gamma:\bR \to M$ is periodic of period $l$, i.e.,\ $\gamma(t + l) = \gamma(t)$ for all $t\in\bR$,
\item Each geodesic is injective on the interval $(0,l)$.
\end{itemize}
A fundamental example of a Zoll manifold is $S^2$ with its round metric. Other known examples are given by the so-called \emph{CROSS} (Compact, Rank-One Symmetric Spaces), consisting of the following five cases:
\begin{enumerate}
\item $S^n$ with the round metric,
\item $\mathbb{RP}^n$ with the induced metric from the round $S^n$,
\item $\bCP$ with the Fubini--Study metric,
\item $\mathbb{HP}^n$ with its canonical metric,
\item $\mathbb{OP}^2$ with its canonical metric.
\end{enumerate}

We record several facts about Zoll manifolds. Since every geodesic on a Zoll manifold is periodic, the manifold is geodesically complete.
Moreover, since all geodesics have the same period, the distance between any two points is at most half the period. Therefore, from the Hopf--Rinow theorem, any Zoll manifold is compact.

Furthermore, Zoll manifolds exhibit some topological rigidity. By the Bott--Samelson theorem, for any Zoll manifold $M$, the fundamental group satisfies $\pi_1(M) = \{0\}$ or $\mathbb{Z}_2$. In the simply connected case, the cohomology ring of $M$ coincides with that of one of the compact rank-one symmetric spaces $S^n$, $\mathbb{CP}^n$, $\mathbb{HP}^n$, or $\mathbb{OP}^2$. If $\pi_1(M) = \mathbb{Z}_2$, then the cohomology ring of $M$ agrees with that of $\mathbb{RP}^n$. We refer to \cite[Chapter 7. D]{B12} for the exposition of these facts. 

Motivated by this theorem, we introduce the following terminology.
\begin{defn}
Let $M$ be a Zoll manifold. If the cohomology ring of $M$ is isomorphic to that of a given CROSS model $M_{\mathrm{model}}$, we say that $M$ is a \emph{Zoll manifold of type $M_{\mathrm{model}}$}. In particular, every Zoll manifold is of one of the five types $S^n$, $\mathbb{RP}^n$, $\mathbb{CP}^n$, $\mathbb{HP}^n$, or $\mathbb{OP}^2$.
\end{defn}

A natural question then arises: for those Zoll manifolds of the five types described, can there exist a Zoll metric different from that of the canonical model? This was first answered in 1903 by Zoll, who constructed a two-parameter family of Zoll metrics on $S^2$. A. Weinstein later extended this line of reasoning to show that such \emph{exotic} Zoll metrics can exist on $S^n$ for $n\ge 2$. L. Green proved that there is no exotic Zoll metric on $\bR\bP^2$. R. Michel proved that,  at least infinitesimally,  no such exotic Zoll metric can arise on $\mathbb{RP}^n$. Very little is currently known about the other three types. We refer again to the classical book \cite{B12} for extensive study on related topics. 

In that context, one naturally wonders under what additional conditions a Zoll metric might coincide with its canonical model. One of the most intriguing geometric features of a Zoll manifold lies in its tangent bundle. Since all geodesics of a Zoll manifold $M$ share the same period, the geodesic flow on $TM$ is itself periodic. Also recall that on any Riemannian manifold $(M,g)$, for the projection $\tau:TM\to M$ and its differential $d\tau:TTM \to TM$, one may define the Liouville $1$-form
\[
\alpha_{(p,v)}(-) \;=\; g_p\bigl(d\tau(-),v\bigr),
\]
and hence the symplectic form $\omega=d\alpha$ on $TM$. Since $\omega$ remains invariant under the geodesic flow, the tangent bundle of a Zoll manifold is a symplectic manifold equipped with an equivariant $S^1$-action. It is thus natural to ask if imposing an additional geometric condition, compatible with this symplectic form---such as, a K\"ahler structure---might force a Zoll metric to be the canonical one. 


In fact, the tangent bundle of each CROSS model has a K\"ahler structure which is compatible with the canonical symplectic structure. According to Patrizio and Wong \cite{PW91}, for each CROSS model $M$ there exists a unique Stein manifold structure on the tangent bundle $TM$ such that for the norm function $u(v) = \|v\|$, $u^2$ is a strictly plurisubharmonic exhaustion function, and $u$ satisfies the homogeneous complex Monge--Ampere equation (HCMA). More precisely, if $m = \dim_\bR M$, then $u$
satisfies
\[
(\sqrt{-1}\,\partial\bar\partial u)^m = 0,
\]
or equivalently,
\[
\det\!\left(
\frac{\partial^2 u}{\partial z_k\,\partial \bar z_l}
\right)_{k,l=1}^{m} = 0
\]
on $TM \setminus M$.


Let $M$ be a real $m$-dimensional Riemannian manifold whose tangent bundle $TM$ admits an integrable complex structure. If the norm function $u$ satisfies the above HCMA and if $u^2$ is a strictly plurisubharmonic exhaustion function on $TM$, then $(TM,M,u^2)$ is called a \emph{Monge--Ampere model}, and the tangent bundle $TM$ is said to admit an \emph{entire Grauert tube}. In this paper, our main result is as follows.



\medskip
\noindent\textbf{Main Theorem.}
\textit{Let $(M,g)$ be a Zoll manifold of type $\bCP$. Assume that $(M, g)$ admits an entire Grauert tube. Then, $M$ is isometric to $\bCP$ with its Fubini--Study metric, up to constant multiplication.}
\medskip
 
Such type of uniqueness theorem for Zoll manifolds of type $S^n$ with an entire Grauert tube was proved by Burns and Leung in \cite{BL18}, and they asked similar questions for Zoll manifolds of other types. The above result answers their question for the case of type $\bCP$. Our proof follows a similar line of argument from \cite{BL18}, and there are some new ingredients. We first use the result of \cite{BL18} to get an analytic compactification $\oX$ of the tangent bundle of $M$ as a complex manifold, and consider the added point set $D$ as an ample divisor. After that, we will calculate the cohomology of $\oX$ and $D$ via algebraic topological techniques, and consider the generators of their second cohomology classes. 
We can show that the complex manifold $\oX$ is a $2n$-dimensional Fano manifold (see Lemma \ref{lem-Fano}) with Fano index $n+1$ (Lemma \ref{lem-ind}). 
In the case considered in \cite{BL18}, this step is relatively easy because the second Betti number of the compactified space is equal to 1, while in our case $b_2(\oX)=2$ and more arguments are needed (see Section \ref{sec-cocycle}). Moreover we will calculate the Fano index of $\oX$ with the help of the Bott-Samelson Theorem and by relating the Morse index of closed geodesics to the Fano index of $\oX$. Finally, by a classification theorem of Wiśniewski \cite{W90}, we can show that $\oX$ is biholomorphic to the standard model $\bC\bP^n\times \bC\bP^n$ (\cite{PW91}). Some algebro-geometric argument allows us to establish that the biholomorphism preserves key structures in $\oX$, which shows that the metric on $M$ should be isometric to the Fubini--Study metric, up to constant multiplication.

We remark that our result does not rule out the existence of exotic Zoll metrics. Indeed,  even though Burns and Leung proved a uniqueness result for Zoll manifolds of type $S^n$ under the entire Grauert tube assumption, the sphere $S^n$ itself admits many non-trivial Zoll structures. On the other hand, our result has the following 
implication from the viewpoint of Stein geometry. Patrizio and Wong~\cite{PW91} proved that a Monge--Ampere model (see section \ref{sec-MAmodel} for an explanation) with center isometric to a compact rank-one symmetric space is unique up to biholomorphic isometry. Our main theorem can be viewed as a strengthening of their result when the center is only assumed to be a Zoll manifold of type $\bCP$, as reflected in the following reformulation.

\medskip
\noindent\textbf{Main Theorem (Rephrased).}
\textit{Let $Y$ be a complex manifold of dimension $2n$ and $u: Y\rightarrow [0, \infty)$ be a plurisubharmonic exhaustion function that form a Monge-Amp\`{e}re model with center $M=\{u=0\}$. Assume that $M$ is a Zoll manifold of type $\bCP$. Then $(Y,M,u^2)$ is unique up to biholomorphic isometry. More precisely, it is isomorphic to $\left((\mathbb{CP}^n\times \mathbb{CP}^n)\setminus D, \mathbb{CP}^n, u_0^2\right)$ where $D=\{\sum_\alpha Z_\alpha W_\alpha=0\}$ (see \eqref{eq-Dmodel}) and $u_0$ is given in \eqref{eq-u0}.} 


\vskip 3mm 

\noindent
\textbf{Acknowledgments.} The authors are grateful to anonymous reviewers for their valuable comments, which greatly improved the manuscript. The authors are partially supported by NSF (Grant No. DMS-2305296).

\section{Complexification and Compactification of Tangent Bundle}\label{sec-MAmodel}
\subsection{Homogeneous Complex Monge--Ampere Equation on the Tangent Bundle}\label{sec-GT}

To investigate the condition under which there is a natural complexification of the entire tangent bundle $TM$ of a manifold $M$ (i.e., an entire Grauert tube), we first need to recall its definition. The basic references are \cite{GS91, LS91}.

The tangent bundle $TM$ of a Riemannian manifold $(M,g)$ has a foliation structure, often referred to as the \emph{Riemann foliation}. For a real number $\tau$, define the dilation operator
\[
N_\tau : TM \,\to\, TM, 
\quad N_\tau(v) \;=\; \tau\,v.
\]
For a geodesic $\gamma$, define the immersion $\phi_\gamma:\bC\rightarrow TM$ by
\[
\phi_\gamma(\sigma+i\tau)\;=\;N_\tau\bigl(\gamma'(\sigma)\bigr).
\]
We remark that we will also use $\phi_\gamma$ to denote the corresponding leaf, i.e., its image in $TM$. If two geodesics $\gamma_1$ and $\gamma_2$ have different images, then under this construction, the images $\phi_{\gamma_1}(\bC\backslash\bR)$ and $\phi_{\gamma_2}(\bC\backslash\bR)$ do not intersect. Consequently, these immersions define a foliation on $TM\backslash M$.

If $T^rM:=\{v \in TM : ||v||_g <r\}$ has a complex structure $J$ such that the map $\phi_{\gamma}$ from $\{z\in\bC:|\Im z| < r\}$ to $T^rM$ is holomorphic, then we call $J$ an \emph{adapted complex structure}, and $T^rM$ a \emph{Grauert tube}. We also set $T^\infty M := TM$. If the above structure exists on $T^\infty M$, then $M$ is said to admit an \emph{entire Grauert tube}.

By Theorem~4.2 of \cite{LS91}, an adapted complex structure is unique whenever it exists. Moreover, Theorem~2.2 of \cite{S91} states that any compact real-analytic Riemannian manifold admits an $\epsilon>0$ such that $T^\epsilon M$ carries an adapted complex structure.

A Grauert tube $T^rM$ has a Stein manifold structure. Define the energy function $E:TM\rightarrow \bR$ by
\[
E(v) \;=\;\frac{1}{2}\|v\|^2.
\]
For $u:=\sqrt{2E}$, $u^2$ is a strictly plurisubharmonic exhaustion function that satisfies the following homogeneous complex Monge--Ampere equation (HCMA)
\[
(i\partial\bar{\partial}u)^m=0,
\]
or equivalently,
\[
\det\!\left(
\frac{\partial^2 u}{\partial z_k\,\partial \bar z_l}
\right)_{k,l=1}^{m} = 0
\]
on $T^rM\backslash M$. The triple $(T^rM,M,u^2)$ is called a \emph{Monge--Ampere model} with center $M$.

Note that $\partial\bar\partial u^2=2u~\partial\bar\partial u+2~\partial u\wedge\bar\partial u$ implies $\operatorname{rank}(\partial\bar\partial u^2)\le\operatorname{rank}(2u~\partial\bar\partial u)+\operatorname{rank}(2~\partial u\wedge\bar\partial u)$. Since $u^2$ is strictly plurisubharmonic, $\operatorname{rank}(\partial\bar{\partial} u^2)=m$. Moreover, $\operatorname{rank}(\partial u\wedge\bar{\partial} u)=1$. Therefore, $\operatorname{rank}(\partial\bar{\partial} u)\ge m-1$, and $(\partial\bar\partial u)^m=0$ forces $\operatorname{rank}(\partial\bar{\partial} u)=m-1$. The kernel of $\partial\bar{\partial} u$ is a rank 1 subbundle of the tangent bundle of $T^rM\backslash M$, implying that $T^rM\backslash M$ admits a foliation by Riemann surfaces \cite{BK77}, which is called the \emph{Monge--Ampere foliation}. According to Theorem 3.1 in \cite{LS91}, on Grauert tubes, the two foliations described above coincide, and each leaf is a complex submanifold under the complex structure $J$.

An adapted complex structure can be concretely constructed from the metric $g$ via \textit{parallel vector fields} on each leaf. Before defining it, recall that for $v\in TM$, a vector $\zeta \in T_vTM$ can be written in terms of its horizontal part $h\in TM$ and vertical part $w\in TM$. More precisely, if  $V:\bR\rightarrow TM$ is a vector field along a path in $M$ such that $V(0)=v$, $(\pi V)'(0)=h$, and $\frac{DV}{dt}(0)=w$, we can express $\zeta=V'(0)\in T_vTM$ as $\zeta=(h,w)$.

For a unit speed geodesic $\gamma:\bR\rightarrow M$, consider the leaf $\phi_\gamma:\bC\rightarrow TM$ of the Riemann foliation. A vector field $V:\bC\rightarrow TTM$ on this leaf is called a \emph{parallel vector field} if it is invariant under the geodesic flow $\Phi_t$ and the dilation $N_\tau$ on $TM$. More concretely, let $\Gamma:(-\epsilon,\epsilon)\rightarrow TM$ be a defining curve of a vector $V(i)=\zeta=(h,w)\in T_vTM$, i.e., $\Gamma'(0)=\zeta$, $(\pi\circ\Gamma)'(0)=h$, $\frac{D}{dt}\Gamma(0)=w$. Then $V(\sigma+i\tau)$ should be the same as $(\Phi_\sigma\circ N_\tau\circ\Gamma)'(0)$. Therefore, a parallel vector field can be determined by its value at $i\in\bC$. Let $\iota:\bR\rightarrow TM$ be the Jacobi field along $\gamma$ with 
$\iota(0)=h$ and $\iota'(0)=w$. 
By definition, in horizontal/vertical notation we then have the equality: 
\[
V(\sigma + i\tau) = \bigl(\iota(\sigma),\tau\,\tfrac{D\iota}{dt}(\sigma)\bigr). 
\]

Now suppose $J:TT^rM\rightarrow TT^rM$ is a well-defined adapted complex structure on $T^rM$, meaning that each leaf of the Riemann foliation is a complex submanifold. By \cite{S91} and \cite{LS91}, such a complex structure $J$ is uniquely determined by the metric $g$. Concretely, let $p\in M$ and a unit vector $v\in T_pM$. Consider the geodesic $\gamma$ with $\gamma(0)=p$ and $\gamma'(0)=v$. Take an orthonormal basis $v_1,\dots,v_{m-1}, v_m=v$ of $T_pM$. Define parallel vector fields $\zeta_i,\eta_i$ by
\[
\zeta_i(i)=(v_i,0), 
\quad \eta_i(i)=(0,v_i)
\]
Then $\zeta_i$ and $\eta_i$ give rise to holomorphic sections of $T^{1,0}(TM)$ on the leaf $\gamma(\bC\backslash\bR)$: $\zeta_i^{1,0}=\zeta_i-\sqrt{-1}J\zeta_i$ and $\eta_i^{1,0}=\eta_i-\sqrt{-1}J\eta_i$. The vector fields $\{\zeta^{1,0}_i; i=1,\dots, m\}$ and $\{\eta^{1,0}_i; i=1,\dots, m\}$ are linearly independent, thus forming a holomorphic basis. For holomorphic vector tuples 
\[\zeta^{1,0}=[\zeta_1^{1,0},\dots,\zeta_{m-1}^{1,0}]^T,\quad \eta^{1,0}=[\eta_1^{1,0},\dots,\eta_{m-1}^{1,0}]^T,\]
we have a holomorphic matrix $\Psi:\bC\backslash\bR\rightarrow M_{m-1}(\bC)$ such that
\[
\eta^{1,0}=\Psi\zeta^{1,0}
\]
On the real axis (except for some discrete set $S\subset\bR$), $\zeta_i$ and $\eta_i$ also do not vanish identically, because they are non-trivial Jacobi fields. For real vectors
\[
\zeta=[\,\zeta_1,\dots,\zeta_{m-1}]^T, 
\quad \eta=[\,\eta_1,\dots,\eta_{m-1}]^T
\]
we have a real matrix $\psi:\bR\backslash S\rightarrow M_{m-1}(\bR)$ such that
\[
\eta \;=\; \psi\,\zeta.
\]
On the entire complex plane $\bC$, we can write
\[
\eta \;=\; \mathrm{Re}\,\Psi\,\zeta \;+\;\mathrm{Im}\,\Psi\,J\zeta,
\]
which implies
\[
J\,\zeta 
\;=\;
(\mathrm{Im}\,\Psi)^{-1}\bigl(\,\eta \;-\;\mathrm{Re}\,\Psi\,\zeta\bigr).
\]
Note that $\psi$ is defined even without the existence of a Grauert tube. Therefore, the existence of a Grauert tube $T^rM$ is essentially equivalent to the condition that the analytic continuation $\Psi$ must be well-defined and $\mathrm{Im} \Psi$ invertible. Indeed, for a compact real-analytic manifold $M$, there exists some $\epsilon>0$ such that $\psi$ can be analytically continued on $T^\epsilon M$.

\subsection{Compactification of the Tangent Bundle of Zoll Manifolds}\label{structure}

Let $M$ be an $m$-dimensional Zoll manifold with period $l$. It is known that the tangent bundle $TM$ admits a compactification
\[
X = TM \cup D
\]
by adjoining the space $D$ of oriented images of geodesics. Historically, $D$ is called the \emph{manifold of oriented geodesics}; see Besse \cite[Chapter 2]{B12}. In the case where $TM$ is an entire Grauert tube, Burns and Leung \cite{BL18} proved that the resulting compactification $X$ is a projective variety. For later use, we first recall their compactification/algebraization result. 

We begin by noting that the set $D$ of oriented geodesics can be identified with the quotient of the unit tangent bundle $UM$ by the equivalence relation induced by the geodesic flow. By the Zoll condition, all geodesic flows are periodic, and hence $D$ is identified with $UM/S^1$ as a set. Since the compact Lie group $S^1$ acts smoothly, freely, and properly on the unit tangent bundle $UM$, the quotient manifold theorem (cf.~\cite{L03}, Theorem~21.10) implies that $D$ inherits a smooth manifold structure as $UM/S^1$. A more detailed description of the chart structure of $D$ will be given later.

We now explain the intuition behind the compactification $X = TM \cup D$. For each leaf of the Riemann foliation (Section \ref{sec-GT}), since the geodesic flow is periodic with period $l$, each leaf $\phi_\gamma$ is biholomorphic to a cylinder without real part,
\[
\phi_\gamma \colon \mathbb{C}/l\mathbb{Z} \setminus \mathbb{R}/l\mathbb{Z} \to TM.
\]
For $\sigma\in\bR$, $\phi_\gamma(\sigma)$ is naturally embedded into the zero section of the tangent bundle, namely
\[
\phi_\gamma(\sigma) = N_0(\gamma'(\sigma)).
\]
Thus, each $\phi_\gamma$ may be regarded as a long cylinder $\mathbb{C}/l\mathbb{Z}$. By definition, for each leaf $\phi_\gamma$, the leaf corresponding to the orientation-reversed geodesic $-\gamma$ is given by a coordinate change induced by the reflection about the origin on $\phi_\gamma$. Consequently, a single leaf $\phi_\gamma$ corresponds bijectively to a pair of elements $(\gamma,-\gamma)$ in $D$. This suggests a compactification in which each cylinder is completed to a $\mathbb{CP}^1$ by adjoining the limits
\[
\gamma = \lim_{\tau \to +\infty} \phi_\gamma(\sigma + i\tau),
\qquad
-\gamma = \lim_{\tau \to -\infty} \phi_\gamma(\sigma + i\tau).
\]
To avoid an abuse of notation, whenever a geodesic $\gamma \in D$ is regarded as a point of $D$, we denote it by $0_\gamma$. Therefore, $\lim_{\tau \to +\infty} \phi_\gamma(\sigma + i\tau)=0_\gamma$ and $\lim_{\tau \to -\infty} \phi_\gamma(\sigma + i\tau)=0_{-\gamma}$. We call this extended leaf the \textit{compactified leaf} $\phi_\gamma$.

This construction can be made rigorous in the following two ways, yielding a smooth structure on the set $X = TM \cup D$.

The first approach uses the symplectic cut construction from \cite{Ler95} (see also \cite{A07}). Observe that $TM$ is diffeomorphic to its open ball bundle over $M$.
Since $D\cong UM/S^1$ where $UM$ is the unit tangent bundle, $X$ as a topological space may be obtained from the closed unit ball bundle $T^1M=\{v\in TM; \|v\|\le 1\}$ by taking the free quotient of its boundary $UM=\partial(T^1 M)$ by the $S^1$-action of the periodic geodesic flow.
On the other hand, equipping $TM\cong T^*M$ with the standard symplectic structure $\omega$, 
we can apply the symplectic cut construction to the Hamiltonian function $H:=\|v\|^2: TM\rightarrow \bR_{\ge 0}$ at level 1 following \cite[1.1: Basic construction]{Ler95}. More precisely, this means that we consider the product symplectic manifold $(TM\times \bC, \omega\oplus \sqrt{-1}dz\wedge d\bar{z})$ equipped with the $S^1$-action:
\begin{equation*}
    e^{i\theta}\circ (v, z)=(\Phi_{l\theta/(2\pi)}(v), e^{-i\theta} z)
\end{equation*}
where $\Phi_t: TM\rightarrow TM$ is the geodesic flow. We can then construct $X$ as a smooth manifold by taking the symplectic quotient $X=\mu^{-1}(1)/S^1$ for the momentum map $\mu:=H+|z|^2: TM\times \bC\rightarrow \bR_{\ge 0}$. This yields a well-defined smooth symplectic structure on $X$ and also realizes $D=H^{-1}(1)/S^1\cong (\mu^{-1}(1)\cap (TM\times \{0\}))/S^1$ as a codimension two symplectic submanifold of $X$. Moreover, there is a natural smooth projection $\pi: X\setminus M\rightarrow D$ induced by the natural map $\bar{p}_1: (TM\setminus M)\times \bC\rightarrow D=UM/S^1$ given by  $\bar{p}_1(v, z)=\left[\frac{v}{\|v\|}\right]$, which realizes $\pi$ as a smooth open disc bundle over $D$. 

We can also give an explicit description by using Riemannian geometry. We will construct local charts of $X$ by extending local charts of $D$ as follows.

\begin{defn}\label{chart}
Let $D$ be the set of oriented images of geodesics in $M$. Let $\gamma$ be a unit-speed geodesic in $M$ with $\gamma(0)=p$ and $\gamma'(0)=v$. An element of $D$ corresponding to $\gamma$ is denoted $0_\gamma$.
Consider the normal vector space $N_vM=\gamma'(0)^\perp \subset T_pM$.
For $u,w\in N_v M$, define a geodesic $\gamma(u,w)$ by
\[
\gamma(u,w)(s):=\exp_{\exp_p u}\!\Bigl( s\, P_u\!\Bigl(\frac{v+w}{\sqrt{1+\|w\|^2}}\Bigr)\Bigr),
\]
where $P_u$ denotes parallel transport along the curve $t\mapsto \exp_p(tu)$ from $p$ to $\exp_p u$. For a sufficiently small open neighborhood $U\subset N_v M$ of $0\in N_v M$, we define a smooth local parametrization $\varphi'_\gamma$ of $D$ around $0_\gamma$ by
\[
\varphi'_\gamma \colon U\times U \subset N_vM\times N_vM \cong \bR^{2m-2} \longrightarrow D,
\qquad
\varphi'_\gamma(u,w):=0_{\gamma(u,w)}.
\]
Let $\Delta:=\{z\in \bC \mid |z|<1\}$. Define $\varphi_\gamma \colon U\times U\times \Delta \longrightarrow X$ by
\[
\varphi_\gamma\!\bigl(u,w,re^{i\theta}\bigr)=
\begin{cases}
\dfrac{l}{4\pi}\cosh^{-1}\left(\dfrac{1}{r^2}\right)\,\gamma(u,w)'\!\left(\dfrac{l}{2\pi}\theta\right), & r\in (0,1),\\[0.4em]
0_{\gamma(u,w)}\in D, & r=0.
\end{cases}
\]
Finally, define $\pi \colon X\setminus M \to D$ by
\[
\pi\bigl(a\,\gamma'(b)\bigr)=0_\gamma \quad \text{for any } a,b>0,
\qquad\text{and}\qquad
\pi(0_\gamma)=0_\gamma.
\]
\end{defn}

The map $\varphi'_\gamma$ is the standard smooth local parametrization of $D$ (see, for instance, Besse~\cite[Chapter 2, 2.14]{B12}). Using this, one proves the following.

\begin{prop}\label{bundle}
Let $M$ be a Zoll manifold, and consider a set $X := TM \cup D$ and $\varphi_\gamma$, $\pi$ defined in Definition \ref{chart}.
\begin{enumerate}
\item The set $X$ admits a unique smooth structure extending the given smooth structure on $TM$ such that each map $\varphi_\gamma$ is a smooth local parametrization.
\item $(X\setminus M, D, \pi, \Delta)$ defines the structure of an open disc bundle. In particular, $X\setminus M$ deformation retracts onto $D$.
\end{enumerate}
\end{prop}

\begin{proof}
Let $\alpha,\beta$ be unit-speed geodesics in $M$ with $\alpha(0)=p$, $\alpha'(0)=v$, and $\beta(0)=p'$, $\beta'(0)=v'$. Consider the maps $\varphi'_\alpha,~\varphi_\alpha,~\varphi'_\beta,~\varphi_\beta$ constructed as in Definition~\ref{chart}. By definition, on the overlap we have
\[
\varphi_\beta^{-1}\circ \varphi_\alpha\!\bigl(u,w,re^{i\theta}\bigr)
=
\Bigl((\varphi_\beta')^{-1}\circ \varphi_\alpha'(u,w),\ re^{i(\theta+t)}\Bigr),
\]
where $t=t(u,w)$ is a function such that for $(\varphi_\beta')^{-1}\circ \varphi_\alpha'(u,w)=(u',w')$,
\[
\alpha(u,w)(0)=\beta(u',w')\left(\dfrac{l}{2\pi}t\right). 
\]
Since each ingredient in this defining relation depends smoothly on $(u,w)$, it follows that $t(u,w)$ is a smooth function.

Writing $z=re^{i\theta}$, we have $re^{i(\theta+t)}=z\,e^{it}$, hence
\[
\varphi_\beta^{-1}\circ \varphi_\alpha(u,w,z)
=
\Bigl((\varphi_\beta')^{-1}\circ \varphi_\alpha'(u,w),\ z\,e^{i t(u,w)}\Bigr),
\]
which is smooth (in particular, smooth at $z=0$). Moreover, the map $\varphi_\gamma$ restricts to a smooth map
\[
\varphi_\gamma\big|_{U\times U\times (\Delta\setminus\{0\})}\colon U\times U\times (\Delta\setminus\{0\})\longrightarrow TM
\]
with respect to the standard smooth structure on $TM$. Hence, for each point of $M\subset TM$, we may take the usual local charts of $TM$ around that point; together with the local parametrizations $\varphi_\gamma$, these charts cover $X$, and all transition maps between them are smooth. Therefore, they determine a unique smooth manifold structure on $X$ by the smooth manifold chart lemma (cf.~\cite[Lemma~1.35]{L03}).

Finally, since $\pi\circ\varphi_\gamma(u,w,z)=0_{\gamma(u,w)}=\varphi'_\gamma(u,w)$, the following diagram commutes:
\[
\begin{tikzcd}[column sep=huge]
\varphi_\gamma(U\times U\times \Delta)
  \arrow[r, "(\varphi'_\gamma\times \mathrm{id}_{\Delta})\circ \varphi_\gamma^{-1}"]
  \arrow[d, "\pi"']
&
\varphi'_\gamma(U\times U)\times \Delta
  \arrow[dl, "\mathrm{proj}_{\varphi'_\gamma(U\times U)}"]
\\
\varphi'_\gamma(U\times U) &
\end{tikzcd}
\]
Therefore, the pairs
\[
\Bigl(\varphi_\gamma(U\times U\times \Delta),\ (\varphi'_\gamma\times \mathrm{id}_{\Delta})\circ \varphi_\gamma^{-1}\Bigr)
\]
form a local trivialization of the disc bundle $(X\setminus M, D, \pi, \Delta)$, proving the second statement. The deformation retraction onto $D$ follows by radial contraction in each disc fiber.
\end{proof}

According to \cite{BL18} and \cite{L14}, we have the following:

\begin{prop}\label{comp}
For a Zoll manifold $M$ with period $l$ admitting an entire Grauert tube, let $\oX=TM\cup D$ be the manifold obtained by compactifying all leaves.
\begin{enumerate}[label=(\roman*)]
\item The complex structure on $TM$ smoothly extends to $X$. (\cite[Lemma 5]{BL18})
\item $\rho=\log\bigl(1+\cosh\bigl(\frac{4\pi\sqrt{2E}}{l}\bigr)\bigr)$ is a K\"ahler potential of $TM$, and the K\"ahler structure derived from this potential extends to $\oX$. (\cite[Proposition III.46]{L14})
\item There exists a Hermitian metric $h$ on $\cO(D)$ such that the curvature form $\Theta_h(\cO(D))$ coincides with the above extended K\"ahler form from $\rho$. In particular, $\cO(D)$ is a positive line bundle. (\cite[Lemma 6]{BL18})
\end{enumerate}
\end{prop}

Therefore, by the Kodaira Embedding Theorem, $D$ is an ample divisor and $\oX$ is a projective manifold. Recall that, in the compactification of the leaf $\phi_\gamma$, we add two points $0_\gamma$ and $0_{-\gamma}$ corresponding to the geodesic $\gamma$ and the same geodesic with opposite orientation $-\gamma$. We note that these points satisfy $0_\gamma=\varphi_\gamma(0,0,0)$ and $0_{-\gamma}=\varphi_{-\gamma}(0,0,0)$, and lie in the smooth images $\varphi_\gamma(0,0,\Delta)$ and $\varphi_{-\gamma}(0,0,\Delta)$. In particular, the compactified leaf $\phi_\gamma$ is a smooth submanifold of $X$. Since the original leaf is a complex submanifold of the entire Grauert tube and the complex structure $J$ extends smoothly to $X$, it follows that the compactified leaf is also a complex submanifold of $X$, biholomorphic to $\mathbb{CP}^1$.

Moreover, we also note the dilation $N_\tau$ with $\tau=-1$, namely $N_{-1}:\oX\rightarrow\oX$. According to \cite{LS91}, $N_{-1}$ is an anti-holomorphic involution on $X$. In particular, $N_{-1}$ fixes exactly the points of $M$ and acts as an orientation-reversing map on $X$, on $D$, and on each leaf $\phi_\gamma$. The existence of this map will be crucial in the proofs that follow.

\subsection{Concrete Description of the Model Case}\label{Model}
Our goal is to show that the compactification $X$ of $TM$, as described in Section \ref{structure}, is biholomorphic to the compactification of the tangent bundle of the canonical model which is $\bCP$ equipped with the Fubini--Study metric. We aim to demonstrate that this biholomorphism preserves the various structures attached to $X$: $M$, $D$, $N_{-1}$, and the exhaustion function $u$. 

To carry out this comparison, we first need to understand precisely how the aforementioned structures are realized in the canonical model of $\bCP$. Patrizio and Wong \cite{PW91} showed that the compactification of $T\bCP$ is biholomorphic to $\tCP$. Building on that, we will give a concrete description of how the leaves, $D$, $u$, $N_{-1}$, and the K\"ahler potential appear in $\tCP$.

Throughout this section, vectors in $\mathbb C^{n+1}$ are viewed as column vectors. First, within $\tCP$, one identifies $\bCP$ via the following totally real embedding:
\[
\begin{aligned}
\iota:\quad & \bCP \;\longrightarrow\;\tCP, \\
& [z]\;\longmapsto\;([\,z\,]\;,\;[\overline{z}\,]).
\end{aligned}
\]
To understand how \(T\bCP\) embeds into \(\tCP\), let us examine the geodesics on \(\bCP\) and their corresponding leaves. Observe that any geodesic \(\gamma\) on \(\bCP\) can be represented using two complex vectors \(z,w\in\bC^{n+1}\) with \(\|z\|=\|w\|=1\) and \(\langle z,w\rangle=\sum_i z_i \bar{w}_i=0\), so that
\[
\gamma(t)\;=\;\bigl[\cos(\frac{t}{2})\,z \;+\;\sin(\frac{t}{2})\,w\bigr],\quad 0\le t\le 2\pi.
\]
Thus,
\[
(\iota\circ\gamma)(t)\;=\;\Bigl(\Bigl[\cos(\frac{t}{2})\,z \;+\;\sin(\frac{t}{2})\,w\Bigr]
,
\Bigl[\cos(\frac{t}{2})\,\overline{z} \;+\;\sin(\frac{t}{2})\,\overline{w}\Bigr]\Bigr).
\]
A natural holomorphic extension of this map is given by
\[
\phi_{\gamma}(\sigma+i\tau)
\;=\;
\Bigl(\Bigl[\cos\frac{\sigma+i\tau}{2}\,z \;+\;\sin\frac{\sigma+i\tau}{2}\,w\Bigr]
,
\Bigl[\cos\frac{\sigma+i\tau}{2}\,\overline{z} \;+\;\sin\frac{\sigma+i\tau}{2}\,\overline{w}\Bigr]\Bigr).
\]
Hence any nonzero vector \(\tau\cdot\gamma'(0)\in T\bCP\) can be represented in \(\tCP\) via
\[
\phi_{\gamma}(i\tau)
\;=\;
\bigl(\bigl[\cosh(\frac{\tau}{2})\,z \;+\; i\,\sinh(\frac{\tau}{2})\,w\bigr]
,
\bigl[\cosh(\frac{\tau}{2})\,\overline{z} \;+\;i\,\sinh(\frac{\tau}{2})\,\overline{w}\bigr]\bigr).
\]
On the compactification, the two points corresponding to the endpoints of each leaf  
are
\[
\lim_{\tau\to\infty}\phi_{\gamma}(i\tau)
\;=\;
0_\gamma\;=\;
([z + i\,w]\;,\;[\overline{z} + i\,\overline{w}]),
\qquad
\lim_{\tau\to -\infty}\phi_{\gamma}(i\tau)
\;=\;
0_{-\gamma}
\;=\;
([z - i\,w]\;,\;[\overline{z} - i\,\overline{w}]).
\]
These added points form the hyperplane section $D\subset \tCP$ (under the Segre embedding) where
\begin{equation}\label{eq-Dmodel}
    D=\left\{([z], [w])\in \tCP; \sum_{0\le \alpha\le n}z_{\alpha}w_{\alpha}=0\right\}.
\end{equation}
Next, let \(N=(n+1)^2-1\) and consider the Segre embedding
\[
\begin{aligned}
S:\quad & \tCP \;\longrightarrow\;\mathbb{CP}^N,\\
& ([z]\;,\;[w]) \;\longmapsto\; [\,z w^T\,].
\end{aligned}
\]
We view homogeneous coordinates $[\zeta]\in \mathbb{CP}^N$ as arranged in the form of an $(n+1)\times(n+1)$ matrix $(\zeta_{ij})_{0\le i,j\le n}$. Define a plurisubharmonic exhaustion $\mathscr{N}:\mathbb{CP}^N \to [1,\infty]$ by
\[
\mathscr{N}(\zeta)
\;=\;
\frac{\sum_{0\le i,j\le n}|\zeta_{ij}|^2}{\bigl|\sum_{0\le\alpha\le n}\zeta_{\alpha\alpha}\bigr|^2},
\]
with the convention that $\mathscr{N}(\zeta)=\infty$ when the denominator vanishes. Indeed, on the open set where $\sum_{\alpha}\zeta_{\alpha\alpha}\neq 0$ we have
\[
\log \mathscr{N}
=
\log\!\Bigl(\sum_{0\le i,j\le n}|\zeta_{ij}|^2\Bigr)
-\log\!\Bigl|\sum_{0\le\alpha\le n}\zeta_{\alpha\alpha}\Bigr|^2,
\]
where the first term is plurisubharmonic and the second term is pluriharmonic. Hence $\log \mathscr{N}$ is plurisubharmonic, and therefore $\mathscr{N}=e^{\log\mathscr{N}}$ is plurisubharmonic. Restricting to \(\tCP \subset \mathbb{CP}^N\), we have
\[
\mathscr{N}(S([z],[w]))
\;=\;
\frac{\sum_{0\leq i,j\leq n}|z_iw_j|^2}{|\sum_{0\leq\alpha\leq n}z_\alpha w_\alpha|^2}
\;=\;
\frac{\sum_{0\le i\le n}|z_i|^2\;\sum_{0\le j\le n}|w_j|^2}{\bigl|\sum_{0\le\alpha\le n}z_{\alpha} w_{\alpha}\bigr|^2}\ge 1.
\]
A straightforward computation shows that
\[
\bigl(\mathscr{N}\circ\phi_{\gamma}\bigr)(i\tau)
\;=\;
\cosh^2(\tau)
\;=\;
\frac{\cosh(2\tau)+1}{2}.
\]
Since on a leaf $\phi_\gamma$, the quantity $|\tau|$ at a point $\phi_\gamma(\sigma+i\tau)$ represents the norm of the corresponding tangent vector, one obtains
\begin{equation}\label{eq-u0}
u_{0} \;:=\; \sqrt{2E}=|\tau|
\;=\;
\tfrac{1}{2}\,\cosh^{-1}\!\bigl(2\,\mathscr{N} - 1\bigr).
\end{equation}

Without loss of generality, for $z=[z_0:\cdots:z_n]$ and $w=[w_0:\cdots:w_n]$, consider the local coordinates with $z_0=w_0=1$. Then the K\"ahler potential $\rho=\log\bigl(1+\cosh\bigl(\frac{4\pi\sqrt{2E}}{l}\bigr)\bigr)$ should be, with $l=2\pi$ and $\sqrt{2E}=u_0$,
$$\log\bigl(2\mathscr{N}(zw^t)\bigr)=\log\bigl(1+\sum_{1\le i\le n}|z_i|^2\bigr)+\log\bigl(1+\sum_{1\le i\le n}|w_i|^2\bigr)-2 \log|\bigl(\sum_{0\le \alpha \le n}z_\alpha w_\alpha\bigr)|+\log 2$$
Since the last two terms are pluriharmonic, they do not contribute to $i\partial\bar{\partial}\rho$, and hence the associated K\"ahler metric coincides with the one induced by the first two terms. These are precisely the Fubini--Study potentials on each factor $\bCP$, respectively. Therefore, the resulting K\"ahler metric is the product of the Fubini--Study metrics on $\tCP$.
\begin{rem}\label{inverse}
The aforementioned embedding of \(T\bCP\) into \(\tCP\) can also be viewed in the reverse direction as follows. Take an arbitrary point \(
([z],[w])\in\tCP\). 

\begin{itemize}
\item If \(\sum_{0 \le \alpha \le n} z_\alpha \,w_\alpha \;\neq\;0,\) then we can choose representatives \(z,w\) such that 
\[
\sum_{0\le \alpha \le n} z_\alpha \,w_\alpha \;=\;1,
\quad
\text{and}
\quad
\|z\| \;=\;\|w\|.
\]
Then \(([z],[w])\) corresponds to a tangent vector at the point 
\(\bigl[\tfrac{z+\overline{w}}{2}\bigr]\in\bCP\). Let $\rho=\frac{\|z+\overline{w}\|}{2}$. Then $2\cosh^{-1}(\rho)$ should be the length of the vector and the direction of the vector should be the direction of the following geodesic.
\[
\Bigl[\cos(\frac{t}{2})\tfrac{z+\overline{w}}{2\rho} \;+\; \sin(\frac{t}{2})\,\tfrac{z-\overline{w}}{2i\sqrt{\rho^2-1}}\Bigr]. 
\]

\item If \(\sum_{0 \le \alpha \le n} z_\alpha \,w_\alpha = 0,\) then pick any representatives \(z,w\) with 
\(\|z\|=\|w\|=\sqrt{2}\). In this case, \(([z],[w])\) can be interpreted as the point in $D$ which represents the geodesic
\[
\Bigl[\cos(\frac{t}{2})\,\tfrac{z+\overline{w}}{2}
\;+\;
\sin(\frac{t}{2})\,\tfrac{z-\overline{w}}{2i}\Bigr]
\]
in \(\bCP\).
\end{itemize}
\end{rem}

\begin{rem}\label{N}\textbf{(Description of $\boldsymbol{N_{-1}}$)}
By definition, for each leaf,
\[
N_{-1}\circ\phi_\gamma(\sigma+i\tau)
\;=\;
\phi_\gamma(\sigma - i\tau).
\]
A direct calculation in this model shows that, in \(\tCP\), 
\[
N_{-1}:\quad
([z],[w]) \;\longmapsto\;([\overline{w}],[\overline{z}]),
\]
and in \(\mathbb{CP}^N\), with the conjugate transpose \(^*\), we have
\[
N_{-1}:\quad
[\zeta]\;\longmapsto\;[\zeta^*].
\]

Notably, all structures we have placed on \(\tCP\) are invariant under \(N_{-1}\). First, the base manifold \(\bCP\) is the fixed locus of \(N_{-1}\) acting on \(\tCP\subset\mathbb{CP}^N\). Furthermore, the involution \(N_{-1}\) sends points of \(\tCP\) to \(\tCP\), points of \(D=\bigl\{\,([z],[w]) : \sum_{0\le \alpha \le n} z_\alpha\,w_\alpha=0\bigr\}\) to \(D\), and each leaf \(L\) to itself respectively.
\end{rem}

\begin{rem}\label{limit}\textbf{(Harmonicity of $\boldsymbol{u}$ on each leaf)}
For any Zoll manifold with entire Grauert tube (not only the model case), consider a geodesic \(\gamma\) corresponding to a point \(0_\gamma\in D\) and its compactified leaf \(\phi_\gamma\). Since this leaf is biholomorphic to \(\mathbb{CP}^1\), for the period $l$ of the geodesic $\gamma$, we may choose a local coordinate \(t = e^{\, \frac{2\pi i}{l}(\sigma + i\tau)}\) centered at \(0_\gamma\). Note that this setup implies \(0_\gamma= \phi_\gamma(i\infty)\), where \(t = 0 = e^{-\infty}\); thus, \(t\) is a holomorphic coordinate (local parameter) near $0_\gamma$ so that $0_\gamma$ is given by $t=0$. On \(\phi_\gamma(\sigma+i\tau)\), the value of \(u\) is \(\lvert\tau\rvert\). Hence,
\[
u \;=\lvert\tau\rvert
\;=\;
\frac{l}{2\pi}\bigl|\log\lvert t\rvert\bigr|,
\]
which shows that \(u\) (restricted to each leaf) is harmonic if $|t|\neq 1$.
\end{rem}

\section{Cohomology of $\oX$ and $D$}

Assume $M$ is a Zoll manifold of type $\bCP$ as before.
First, we compute the cohomology of $\oX$ and $D$, and analyze the behavior of their cocycles. This computation is influenced by \cite{R90} and \cite{Y91}; see also \cite{A07}.

\subsection{Computation of Cohomology}
Before computing the cohomology of $X$, we first remark on the following structural facts.

\begin{rem}\label{TopStructure}
Let $M$ be a Zoll manifold of type $\bCP$, and consider its tangent bundle $TM$, unit vector bundle $UM$,
and the compactification $X := TM \cup D$.
\begin{itemize}
\item The bundle $TM$ deformation retracts onto $M$, $TM \setminus M$ deformation retracts onto $UM$, and by Proposition \ref{bundle}, $X\setminus M$ deformation retracts onto $D$.
\item Recall the following real dimensions:
\[
\dim_\bR M = 2n, \qquad \dim_\bR UM = 4n-1, \qquad \dim_\bR D = 4n-2, \qquad \dim_\bR X = 4n
\]
\end{itemize}
\end{rem}


\begin{prop}\label{cohom}
\noindent
Let $M$ be a Zoll manifold of type $\bCP$, and let $UM$ be the unit vector bundle of $M$. Let $D$ be the space of oriented geodesics of $M$ defined in Section~\ref{structure}, and set $X := TM \cup D$, the compactification of $TM$ as in Proposition~\ref{bundle}. Then the following isomorphisms hold.
\begin{enumerate}
\item 
$H^j(UM,\bZ)\;\cong\;
\begin{cases}
\bZ & j=2k,\;0 \leq k \leq n-1, \text{or}\;j=2k-1,\;n+1\leq k\leq 2n\\
\bZ/(n+1)\bZ & j=2n,\\
0 & \text{otherwise},
\end{cases}$

\item 
$H^j(D,\bZ)\;\cong\;
\begin{cases}
\bZ^{\,\frac{j}{2}+1} & \text{if $j$ is even and } 0 \leq j \leq 2n-2,\\
\bZ^{\,2n-\frac{j}{2}} & \text{if $j$ is even and } 2n \leq j \leq 4n-2,\\
0 & \text{otherwise},
\end{cases}$

\item 
$H^j(\oX,\bZ)\;\cong\;H^j\bigl(\bCP\times\bCP,\bZ\bigr)$ for any $j\in \bZ$. 
\end{enumerate}
\end{prop}

\noindent \textit{Proof.}
Since $M$ is of type $\bCP$, by the Bott-Samelson Theorem, $\pi_1(M)=0$. Therefore, for $n=1$, $M$ is a simply connected $2$-dimensional manifold, namely a sphere, which is resolved in \cite{BL18}. Hence $UM\cong\mathbb{RP}^3$, $D\cong S^2$, $\oX\cong S^2\times S^2$. Now we assume that $n\geq 2$.
First, recall Remark \ref{TopStructure} and observe that $UM$ can be viewed as fiber bundles in two different ways: for the projection $\tau:TM\to M$ and $\pi:X\setminus M \to D$ defined in Definition \ref{chart},
\[
S^{2n-1} \;\hookrightarrow\; UM \;\overset{\tau}{\twoheadrightarrow}\; M
\quad\quad\text{and}\quad\quad
S^1 \;\hookrightarrow\; UM \;\overset{\pi}{\twoheadrightarrow}\; D.
\]
By the long exact sequence of homotopy groups for each fibration, we obtain
\[
\pi_1(S^{2n-1})=0 \;\longrightarrow\; \pi_1(UM) \;\longrightarrow\; \pi_1(M)=0
\]
Hence $\pi_1(UM)=0.$ Also,
\[
\pi_1(UM)=0\;\longrightarrow\;\pi_1(D)\;\longrightarrow\;\pi_0(S^1)\;\longrightarrow\; \pi_0(UM),
\]
which implies $\pi_1(D)=0.$ Since $\pi_1(M)=0$ and $\pi_1(D)=0$, the local coefficient systems arising in the Gysin sequences of the above sphere bundles are trivial. Therefore, the Gysin exact sequences can be applied with integer coefficients.

\noindent\textbf{1. Cohomology of $UM$.} From the Gysin exact sequence associated to the first fibration, we have
\[
H^{j-2n}(M)\;\longrightarrow\;H^j(M)\;\stackrel{\tau^*}{\longrightarrow}\;H^j(UM)\;\longrightarrow\;H^{j-2n+1}(M)\;\longrightarrow\;H^{j+1}(M).
\]
Note that we know that $H^j(M)\cong H^j(\bCP)$ for any $j\in \bZ$. 
Thus, for \(0 \leq j \leq 2n-2\), it follows that 
\[
H^j(UM)\;\cong\; H^j(M)\cong H^j(\bCP),
\]
and for \(2n+1 \leq j \leq 4n-1\),
\[
H^j(UM)\;\cong\; H^{j-2n+1}(M)\cong H^{\,j-2n+1}(\bCP).
\]

Next, consider the case \(j=2n-1\). The part of the Gysin sequence becomes
\[
\underbrace{H^{2n-1}(M)}_{\cong 0}
\;\longrightarrow\;
H^{2n-1}(UM)
\;\longrightarrow\;
H^0(M)\;\xrightarrow{\;\smile e\;}\;H^{2n}(M)
\;\longrightarrow\;
H^{2n}(UM)\;\longrightarrow\;
\underbrace{H^1(M)}_{\cong 0}.
\]
Here $e$ is the Euler class of $TM$ which satisfies $e([M])=\chi(M)=\chi(\bCP)=n+1$.

Consequently,
\[
H^{2n-1}(UM)\;\cong\;\mathrm{Ker}(\smile e)\;\cong\;0,
\]
\[
H^{2n}(UM)\;\cong\;\mathrm{Coker}(\smile e)\;\cong\;\bZ/(n+1)\bZ.
\]
This completes the computation of \(H^j(UM)\).

\noindent\textbf{2. Cohomology of $D$.} We use the Gysin exact sequence for the fibration 
\(\,S^1\hookrightarrow UM \overset{\pi}{\twoheadrightarrow} D.\)
From
\[
H^{2j}(UM)\;\longrightarrow\;H^{2j-1}(D)\;\longrightarrow\;H^{2j+1}(D)\;\longrightarrow\;H^{2j+1}(UM),
\]
and the facts \(H^{-1}(D)\cong 0\) and \(H^{4n-1}(D)\cong 0\), we deduce \(H^{2j+1}(D)=0\) for \(0 \leq j \leq n-1\) by an induction from \(j=0\) to \(j=n-1\) via the last three terms, and for \(j \geq n\) by an induction from \(j=2n-1\) to \(j=n+1\) via the first three terms.

Since \(D\) is a connected real \((4n-2)\)-dimensional manifold, we have \(H^0(D)\cong H^{4n-2}(D)\cong \bZ\). Consider then the Gysin exact sequence
\[
\underbrace{H^{2j-1}(D)}_{\cong 0}
\;\longrightarrow\;
H^{2j-1}(UM)\;\longrightarrow\;
H^{2j-2}(D)\;\longrightarrow\;
H^{2j}(D)\;\stackrel{\pi^*}{\longrightarrow}\;
H^{2j}(UM)\;\longrightarrow\;
\underbrace{H^{2j-1}(D)}_{\cong 0}.
\]
For \(0 \leq j \leq n-1\), we have \(H^{2j}(D)=\bZ^{\,j+1}\) by an inductive argument on the right five terms (from \(j=1\) to \(j=n-1\)). For \(j \geq n\), we have \(H^{2j}(D)=\bZ^{\,2n-j}\) by an inductive argument on the left five terms (from \(j=2n-1\) to \(j=n+1\)). This completes the computation of \(H^j(D)\).

\noindent\textbf{3. Cohomology of $X$.} Since $X=TM\cup(X\backslash M)$, $TM\simeq M$, $X\backslash M \simeq D$, and $TM\backslash M\simeq UM$, we apply the Mayer--Vietoris sequence to compute $H^j(\oX)$.

\smallskip
\noindent\textbf{Case 1. $j\le 2n$.}
Consider the sequence
\[
H^{j-1}(UM)\;\xrightarrow{\;\partial\;}\;H^j(\oX)\;\longrightarrow\;H^j(D)\oplus H^j(M)
\;\xrightarrow{\;\pi^*-\tau^*\;}\;
H^j(UM)\;\xrightarrow{\;\partial\;}\;H^{j+1}(\oX).
\]
The maps $\tau^*$ are surjective, so the boundary maps $\partial$ are zero. Therefore, the middle three terms form a short exact sequence. Hence, if $j$ is odd then $H^j(\oX)=0$. If $j=2k<2n$, then
\[
0\to H^{2k}(\oX)\to \mathbb{Z}^{k+2}\to \mathbb{Z}\to 0,
\]
and therefore $H^{2k}(\oX)\cong \mathbb{Z}^{k+1}$. If $j=2n$, then
\[
0\to H^{2n}(\oX)\to \mathbb{Z}^{n+1}\to \mathbb{Z}/(n+1)\mathbb{Z}\to 0,
\]
and hence $H^{2n}(\oX)\cong \mathbb{Z}^{n+1}$.

\smallskip
\noindent\textbf{Case 2. $j\ge 2n+1$.}
Consider the sequence
\[
H^{j-1}(D)\oplus H^{j-1}(M)\xrightarrow{\;\pi^*-\tau^*\;}H^{j-1}(UM)\xrightarrow{\;\partial\;}
H^{j}(\oX)\longrightarrow H^{j}(D)\oplus H^{j}(M)
\xrightarrow{\;\pi^*-\tau^*\;}H^{j}(UM).
\]
For $j\ge 2n+2$, the group $H^{j-1}(D)\oplus H^{j-1}(M)$ is nonzero only when $j$ is odd, while $H^{j-1}(UM)$ is nonzero only when $j$ is even. Therefore, in this range the map $\pi^*-\tau^*$ is the zero map, and the middle three terms again form a short exact sequence. Hence, if $j\ge 2n+2$ is odd then $H^j(\oX)=0$. If $j=2k \ge 2n+2$ is even, then
\[
0\to \mathbb{Z}\to H^{2k}(\oX)\to \mathbb{Z}^{2n-k}\to 0,
\]
so $H^{2k}(\oX)\cong \mathbb{Z}^{2n-k+1}$.

Finally, for $j=2n+1$, the map $\pi^*-\tau^*$ on the left of the above sequence is surjective, so $\partial=0$. Therefore, $H^{2n+1}(\oX)=0$, which completes the computation.
\qed

\subsection{Behavior of Cocycles}\label{sec-cocycle}
A crucial step in our proof is to show that $\oX$ is a Fano manifold and that $-K_{\oX} = \cO\bigl((n+1)D\bigr)$. To this end, we must first prove that for some integer $r$, 
\[
-K_{\oX} \;=\;\cO(rD)
\quad\text{or equivalently}\quad
c_1(\oX)=r\,[D].
\]
By Proposition~\ref{cohom}, the inclusion \(\iota:D\rightarrow\oX\) induces an isomorphism 
\(\iota^*:H^2(\oX,\bZ)\rightarrow H^2(D,\bZ)\). Hence, it suffices to show that
\[
c_1(D)\;=\;r\,[D]\bigl|_D.
\]
In addition, since \(\cO(D)\bigl|_D \cong \cN_{D/\oX}\) and \(e\bigl(\cN_{D/\oX}\bigr)=c_1\bigl(\cN_{D/\oX}\bigr)\), it follows that
\[
[D]\bigl|_D
\;=\;
c_1\bigl(\cO(D)\bigl|_D\bigr)
\;=\;
c_1\bigl(\cN_{D/\oX}\bigr)
\;=\;
[\omega].
\]
Here, $[\omega]$ is the Euler class in $H^2(D,\bZ)$ induced by the fibration $UM \;\overset{\pi}{\twoheadrightarrow}\; D$. Therefore, we seek explicit expressions for the generators of \(H^2(D,\bZ)\cong \bZ\oplus \bZ\). Observe the following short exact sequence:
\[
0
\;\longrightarrow\;
H^0(D)
\xrightarrow{\;\cup\,\omega\;}
H^2(D)
\xrightarrow{\;\pi^*\;}
H^2(UM)
\;\longrightarrow\;
0.
\]
Hence, for the generator \(y\in H^2(UM)\cong \bZ\), we see that \(H^2(D)\cong \bZ\oplus \bZ\) is generated by the class \([\omega]\) and \(x\in H^2(D)\) such that \(\pi^*(x)=y\) . 

Next, recall the antiholomorphic involution \(N_{-1}:\oX\rightarrow\oX\),
\[
N_{-1}(v)\;=\;-v,\quad
N_{-1}(0_\gamma)\;=\;0_{-\gamma},\quad
N_{-1}(0_{-\gamma})\;=\;0_\gamma.
\]
We note the following properties of \(N_{-1}\):
\begin{enumerate}
\item 
\(N_{-1}\) sends points of \(D\) to points of \(D\) and reverses the orientation on \(D\). Therefore,
\[
N_{-1}^*\,c_1(D)
\;=\;
N_{-1}^*\,c_1(-K_D)
\;=\;
-\,c_1(D).
\]


\item
Recall the correspondence $[\omega]=c_1(\cO(D)|_D)$. Again, since $N_{-1}$ reverses the orientation on $D$, we have
\[
N_{-1}^*[\omega]
\;=\;
N_{-1}^*c_1(\cO(D)|_D)
\;=\;
-c_1(\cO(D)|_D)
\;=\;
-[\omega]
\]

\item 
\(N_{-1}\) fixes the base manifold \(M\). Since \(y\in H^2(UM)\) is induced by a generator \(z\in H^2(M)\) i.e., $\pi^*x=y=\tau^*z$, it follows that \(N_{-1}^*(y)=y\). Thus,
\[
\pi^*\!\bigl(x - N_{-1}^*x\bigr)
\;=\;
y - N_{-1}^*(y)
\;=\;
0.
\]
So for some integer \(d\), we obtain \(x - N_{-1}^*x = d\,\omega\) and the matrix of $N_{-1}^*$ with respect to the basis $\{x, \omega\}$ is given by $
\left(\begin{array}{cc}
    1&-d\\
    0&-1
\end{array}
\right).
$
\end{enumerate}
Consequently, \(c_1(D)\) and \(\omega\) both belong to the \((-1)\)-eigenspace, which is one-dimensional. Since \(\omega\) is one of the generators, there exists an integer \(r\) such that $c_1(D)=r[\omega]$. Therefore, we conclude the following:

\begin{lem}\label{lem-Fano}
\label{fano}
We have the identity:
$-K_{\oX} \;=\;\cO(rD)$  for some integer $r$.
\end{lem}

\section{Evaluation of Fano Index}

One of the powerful features of having an entire Grauert tube is that all compactified leaves of the Riemann foliation become complex submanifolds. This property can be used effectively in computing the Fano index by applying the adjunction formula to a compactified leaf~\(C\). Let us consider a geodesic \(\gamma\) on \(M\), the corresponding leaf \(\phi_\gamma\) of the Riemann foliation, and its compactification \(C\). Denote by \(\cN\) the holomorphic normal bundle of \(C\) in \(\oX\). By the adjunction formula, we have
\[
K_C \;=\; K_{\oX}\bigl|_{C} \;+\;\det \cN.
\]
To exploit this identity, we first prove the following lemma regarding \(\det \cN\). Let $2n$ be the complex dimension of $\oX$.

\begin{lem}\label{normal}
\(\det \cN \cong \cO_{\mathbb{CP}^1}(2n).\)
\end{lem}

\noindent\textit{Proof.}
Fix a unit-speed closed geodesic \(\gamma\). Since \(M\) is of type \(\bCP\), the Bott--Samelson theorem implies that the index \(k\) of \(\gamma\) is \(1\). In other words, there exist unit vectors \(v_1,v_2,\ldots,v_{2n-1}\in T_{\gamma(0)}M\) all orthogonal to \(\gamma'(0)\), and corresponding Jacobi fields \(\eta_i\) on \(\gamma\) such that
\[
\eta_i(0) \;=\; 0,
\quad
\frac{D\eta_i}{dt}(0)\;=\;v_i,
\]
and among these, \(\eta_1\) has exactly one additional first-order vanishing point along \(\gamma\) (besides \(t=0\)), while \(\eta_2,\ldots,\eta_{2n-1}\) vanish to first order only at \(t=0\).

From these Jacobi fields, we obtain parallel vector fields over \(C\). Since each \(v_i\) is orthogonal to \(\gamma'(0)\), the associated parallel vector fields define holomorphic sections \(\eta_i^{1,0}\) of \(\cN\). Consequently, their wedge product \(\bigwedge\eta_i^{1,0}\) gives a holomorphic section of \(\det \cN\). As noted above, \(\eta_1\) vanishes to first order at two points on \(\gamma\), and each \(\eta_2,\ldots,\eta_{2n-1}\) vanishes to first order at exactly one point. Moreover, by the construction of complex structures from \cite{BL18}, $\{\eta_i^{1,0}, i=1,\dots, n\}$ extend over $D\subset X$ to holomorphic sections of $T^{(1,0)}X$ and are linearly independent at points of $D$. In fact, at points of $D$, they form a basis of the tangent space $T^{1,0}D$. 
Thus, the holomorphic section \(\bigwedge \eta_i^{1,0}\) has a total vanishing order \(2n\). Hence, \(\det \cN = \cO_{\mathbb{CP}^1}(2n)\). 
\qed

Using this, we obtain the following result (cf.\ \cite[Lemma~8]{BL18}):

\begin{lem}\label{lem-ind}
We have the identity:
$
-K_{\oX} \;\cong\;\cO\bigl((n+1)D\bigr).
$
In particular, $X$ is a Fano manifold, i.e., $-K_X$ is an ample line bundle. 
\end{lem}

\noindent\textit{Proof.}
By the adjunction formula,
\[
K_C \;=\; K_{\oX}\bigl|_{C} \;+\; \det \cN.
\]
Since \(C\) is biholomorphic to \(\mathbb{CP}^1\), we have \(K_C = \cO(-2)\). From Lemma~\ref{normal}, \(\det \cN = \cO_{\mathbb{CP}^1}(2n)\), implying
\[
K_{\oX}\bigl|_{C} \;\cong\; \cO\bigl(-2n-2\bigr).
\]
Moreover, the divisor \(D\) meets \(C\) transversally at two points \(\{0_\gamma,0_{-\gamma}\}\), so 
\[
\cO(D)\bigl|_{C} \;\cong\;\cO(2).
\]
By Lemma~\ref{fano}, we already know that \(-K_{\oX}\cong\cO(rD)\) for some integer \(r\). Restricting to \(C\), we get
\[
\cO\bigl(-(2n+2)\bigr)
\;\;\cong\;\;
K_{\oX}\bigl|_{C}
\;\;\cong\;\;
-\cO(rD)\bigl|_{C}
\;\;=\;\;\cO(-2r).
\]
Hence \(\cO(-2n-2)\cong\cO(-2r)\), implying \(r=n+1\). Thus the Fano index is \(n+1\). 
\qed

\section{Characterization of the Complex Structure and Metric}

According to Theorem~B of Wi\'sniewski \cite{W90}, any (complex) \(2n\)-dimensional Fano manifold with Fano index \(n+1\) and Picard number \(2\) is biholomorphic to \(\tCP\). By the Kodaira vanishing theorem, the map 
\[
c_1:H^1\bigl(\oX,\cO^*\bigr)\;\longrightarrow\; H^2(\oX,\bZ)
\]
is an isomorphism, and from Proposition~\ref{cohom} we have \(H^2(\oX,\bZ)\cong \bZ\oplus\bZ\). Thus the Picard number of \(\oX\) is \(2\). It follows that:

\begin{prop}\label{algeb}
There is a biholomorphism
\(\oX \;\cong\;\tCP\) and $D$ is identified with a smooth divisor in $\cO_{\oX}(1,1)=p_1^*\cO_{\bCP}(1)\otimes p_2^*\cO_{\bCP}(1)$ where $p_i$, $i=1,2$ are projections to the two factors of $\bCP$. 
\end{prop}



\noindent\textit{Proof of the Main Theorem.}
By Proposition \ref{algeb}, 
we obtain \(\oX \cong \tCP\), together with projections \(p_i:\oX \to \bCP\) \((i=1,2)\) such that 
\[
\cO(D) \;\cong\; p_1^*\cO(1)\,\otimes\,p_2^*\cO(1).
\]
Consider $\phi=\nu\circ p_1\circ N_{-1}$ where $\nu: \bCP\rightarrow \bCP$ is the canonical conjugation involution on $\bCP$. Then $\phi$ is a surjective holomorphic morphism $\oX\rightarrow \bCP$. By Mori's theory, there are only two such holomorphic morphisms $p_1$ and $p_2$ up to isomorphism which correspond to the contraction of two extremal rays. In other words, there exists a $\sigma \in \mathrm{Aut}(\bCP)$ such that either $\phi=\sigma\circ p_2$ or $\phi=\sigma\circ p_1$. 

We consider the first case: $\phi=\nu\circ p_1\circ N_{-1}=\sigma\circ p_2$ with $\sigma\in \mathrm{Aut}(\bCP)$. Replacing $p_2$
by $\sigma\circ p_2$, we can assume that $\nu\circ p_1\circ N_{-1}=p_2$ which implies
$\nu\circ p_2 \circ N_{-1}=p_1$.
Then we have $
p_1\circ N_{-1}([Z],[W])=\nu \circ p_2([Z], [W])=\nu([W])=[\overline{W}]
$
and similarly $p_2\circ N_{-1}([Z],[W])=[\overline{Z}]$. So we get the expression for $N_{-1}$:
\begin{equation*}
    N_{-1}([Z], [W])=([\overline{W}], [\overline{Z}]). 
\end{equation*}
To proceed, we want to exclude the second case. Suppose $\nu\circ p_1\circ N_{-1}=\sigma_1\circ p_1$ for some $\sigma_1\in \mathrm{Aut}(\bCP)$. By the same reasoning as above, we must have $\nu\circ p_2\circ N_{-1}=\sigma_2\circ p_2$ for some $\sigma_2\in \mathrm{Aut}(\bCP)$. So we get:
\begin{equation*}
N_{-1}([Z], [W])=(\overline{\sigma_1([Z])},\; \overline{\sigma_2([W])}).
\end{equation*}
Since $N_{-1}$ is an anti-holomorphic involution,
both $\nu \circ \sigma_i$ are anti-holomorphic involutions of $\bCP$.
The fixed point set $M$ of $N_{-1}$ is equal to 
$K_1\times K_2$ where $K_i$, $i=1,2 $ is the fixed point set of the anti-holomorphic involution $\nu\circ \sigma_i$ which are smooth submanifolds. 
The anti-holomorphic involutions of $\bCP$ are well understood and we know that $K_i$ is either diffeomorphic to $\mathbb{R}\mathbb{P}^n$ or is empty. This would contradict the assumption that $M$ has the same cohomology ring as $\bCP$. 
The latter claim follows from \cite[Lemma 5.29]{SK22}, which states that any anti-holomorphic involution is projectively equivalent to the standard complex conjugation or when $n+1=2m$ is even to the involution given by 
$$
[Z_1,\cdots, Z_{2m}]\mapsto [-\overline{Z_{m+1}}, \cdots, -\overline{Z_{2m}}, \overline{Z_1},\cdots, \overline{Z_m}].
$$
To explain this further, note that since $\mathrm{Aut}(\bCP)=PSL(n+1, \bC)$,  $\sigma=\sigma_A$ for an invertible $(n+1)\times (n+1)$-matrix $A$ and the anti-holomorphic involution $\nu\circ \sigma$ is given by $[Z]\mapsto [\overline{A v}]=[\bar{A}\bar{v}]$. 
Then there exists an invertible $(n+1)\times (n+1)$-matrix $U$ such that either $A=\bar{U}^{-1}{U}$ or $A=\bar{U}^{-1}J{U}$ with $J=\left(\begin{array}{cc}0&-I\\I&0\end{array}\right)$. Note that the latter case happens only when $(n+1)=2m$ is even. 
With this normal form, we see that $[v]$ with $v\in \bC^{n+1}\setminus \{0\}$ is a fixed point of $\nu\circ \sigma_A$ if and only if 
$\overline{Av}=\lambda v$ with $\lambda \neq 0$, which is equivalent to either 
$U v=\bar{\lambda} \overline{Uv}$ or $J U v=\bar{\lambda} \overline{U v}$. We easily see that the fixed point set of $\nu\circ\sigma$ is either diffeomorphic to $\mathbb{R}\mathbb{P}^n$ or is empty. 


From now on, we can assume the equality \(p_2=\overline{\,p_1\circ N_{-1}}\). Assume $D\in H^0(\bCP\times \bCP, \cO(1,1))$ is given by the equation $\sum_{i,j} c_{ij}Z_iW_j=0$. Because $D$ is invariant under the conjugation $N_{-1}: ([Z], [W])\mapsto ([\bar{W}], [\bar{Z}])$, it is easy to see that $C=(c_{ij})$ can be chosen to be a Hermitian matrix $C^*=C$. Therefore, we can diagonalize the matrix $C$ by a unitary matrix $U$: $U C U^*=I$ where $U^*=\bar{U}^t$. By a simultaneous change of variables $Z\mapsto U^t Z$, $W\mapsto \bar{U}^t W$, we can assume $D$ is given by the equation $\sum_i Z_i W_i=0$. 
Hence, under the new coordinate for the Segre embedding \(\oX \to \mathbb{P}\,(H^0\bigl(\oX,\cO(D)\bigr)^*)\), the data $(N_{-1}, D, M)$ coincide with the data of standard model $\bCP$. In particular, note that we now have $M$ is diffeomorphic to $\bCP$.

We now verify that the norm functions $u_0$ and $u_1$, induced respectively by the
canonical metric and by the Zoll metric on $M$, coincide on
$X \cong \tCP$.
Once this is established, the proof is complete, since for the plurisubharmonic
exhaustion $u^2$ of the Monge--Ampere model $(X,M,u^2)$,
the restriction of the K\"{a}hler metric $2 (\sqrt{-1}\partial\bar{\partial}(u^2))\circ J$ to $M$ coincides with
the original Riemannian metric $g$ on $M$ (see \cite[Section~3]{PW91} or \cite[(1.4)]{GS91}).

To verify that $u_0=u_1$, 
the idea is that we restrict $u_1$ to any leaf associated to $u_0$ and argue that $u_1\le u_0$ by using the subharmonicity and the asymptotic behavior of $u_1-u_0$ on the leaf. Since the leaves associated to $u_1$ cover $X\setminus D$, we conclude that $u_1\le u_0$. Reversing the role of $u_0$ and $u_1$, we also get $u_0\le u_1$. We now write down the details. 

Let $(i,j)\in\{(0,1),(1,0)\}$.
Consider a geodesic $\gamma$ of $M$ with respect to the metric associated to $u_i$,
and the corresponding compactified leaf $\phi_\gamma$.
Following the notation of Remark~\ref{limit}, we introduce the coordinate
\[
t = e^{\frac{2\pi i}{l}(\sigma+i\tau)}
\]
on $\phi_\gamma$.

\begin{enumerate}
\item
For the function $u_i$, as observed in Remark~\ref{limit}, we have
\[
u_i|_{\phi_\gamma} = |\tau|
= \frac{l}{2\pi}\,|\log|t||.
\]

\item
For the function $u_j$, the function
\[
\rho = \log\!\left(1+\cosh\!\left(\frac{4\pi u_j}{l}\right)\right)
\]
is one of the possible K\"ahler potentials on $X$ (Proposition~\ref{comp} (ii)).
Accordingly, there exists a Hermitian metric $\|\cdot\|_h$ on $\mathcal O(D)$ such that for the defining section $s_D$ of the divisor $D\subset X$ 
\[
\rho=-\log \|s_D\|_h^2, \quad 
\text{ and } \quad
\Theta_h(\mathcal O(D)) = i\,\partial\bar\partial \rho. 
\]
Since the holomorphic leaf $\phi_\gamma$ meets $D$ transversely and $D|_ {\phi_\gamma}=\{t=0\}$ near $0_\gamma\in \phi_\gamma\cong \mathbb{CP}^1$, we have
\[
\|s_D\|^2_h\big|_{\phi_\gamma} = c(t)\,|t|^2
\]
for some smooth function $c$ with $c(0)\neq 0$.
Therefore, in a neighborhood of $t=0$,
\[
u_j = -\frac{l}{2\pi}\log|t| + O(1)
= u_i + O(1),\qquad u_j-u_i=O(1). 
\]
The same conclusion holds at $t=\infty$ by reversing the orientation of the
geodesic $\gamma$.
\end{enumerate}

Since $u_j$ is plurisubharmonic on $X$, it is subharmonic on the holomorphic leaf $\phi_\gamma$.
On the other hand, as noted in Remark~\ref{limit}, the function $u_i$ is harmonic
on $\{|t|\neq 1\}$.
Consequently, $u_j-u_i$ is subharmonic on $\{|t|\neq 1\}$.
Moreover, since both $u_i$ and $u_j$ vanish on $M$, we have $u_i(t)=u_j(t)=0$ for $|t|=1$.

Define $M(r):=\max\{(u_j-u_i)(z)\,;\ |z|=r\}$. By the Hadamard three-circle theorem
(\cite[Chapter~2, Theorem~28]{PW67}), for $0<r_1<r<r_2$,
\[
M(r)\le
\frac{\log r_2-\log r}{\log r_2-\log r_1}\,M(r_1)
+\frac{\log r-\log r_1}{\log r_2-\log r_1}\,M(r_2).
\]
Since $M(1)=0$ and $M(r)$ is bounded as $r\to0$ and $r\to\infty$, taking the limits
$(r_1,r_2)\to(0,1)$ and $(r_1,r_2)\to(1,\infty)$ yields $M(r)\le 0$ for all $r$.
It follows that $u_j \le u_i$ on $\phi_\gamma$. Since $\phi_\gamma$ was arbitrary, we conclude that $u_j \le u_i$ globally.
Finally, as $(i,j)\in\{(0,1),(1,0)\}$ was arbitrary, this implies $u_0 = u_1$. \qed

\vskip 3mm
\noindent
Department of Mathematics, Rutgers University, Piscataway, NJ 08854-8019.

\noindent
{\it Email:} chi.li@rutgers.edu\\
{\it Email:} ks1951@rutgers.edu


\begin{thebibliography}{999999}

\bibitem{A07}Audin, M. Lagrangian skeletons, periodic geodesic flows and symplectic cuttings. {\em Manuscripta Mathematica}. \textbf{124}, 533-550 (2007)

\bibitem{BK77}Bedford, E. and Kalka, A. Foliations and complex Monge--Ampere equations. {\em Communications On Pure And Applied Mathematics}. \textbf{30}, 543-571 (1977)

\bibitem{B12}Besse, A. Manifolds all of whose geodesics are closed. (Springer Science \& Business Media,2012)

\bibitem{BL18}Burns Jr, D. and Leung, K. The complex Monge--Ampere equation, Zoll metrics and algebraization. {\em Mathematische Annalen}. \textbf{371}, 1-40 (2018)

\bibitem{GS91}
Guillemin, V. and Stenzel, M. Grauert tubes and the homogeneous Monge--Ampere equation. {\em J.
Differential Geom.} \textbf{34} (1991), no. 2, 561–570.

\bibitem{L03}Lee, J.~M. \emph{Introduction to Smooth Manifolds}, Graduate Texts in Mathematics, Vol.~218, Springer, New York, 2003.


\bibitem{LS91}Lempert, L. and Szőke, R. Global solutions of the homogeneous complex Monge--Ampere equation and complex structures on the tangent bundle of Riemannian manifolds. {\em Mathematische Annalen}. \textbf{290}, 689-712 (1991)

\bibitem{Ler95}
Lerman, E. Symplectic cuts. {\em Math. Res. Lett.}, \textbf{2}(3):247–258 (1995). 


\bibitem{L14}Leung, K. Complex Geometric Invariants Associated to Zoll Manifolds. Thesis, University of Michigan (2014)

\bibitem{PW91}Patrizio, G. and Wong, P. Stein manifolds with compact symmetric center. {\em Mathematische Annalen}. \textbf{289}, 355-382 (1991)

\bibitem{PW67}
Protter, M.H. and Weinberger, H.F. Maximum Principles in Differential Equations, Prentice–Hall, Englewood Cliffs, NJ, 1967.


\bibitem{R90}Reznikov, A. On the volume of certain manifolds with closed geodesics. {\em Journal Of Soviet Mathematics}. \textbf{48} pp. 83-86 (1990)


\bibitem{SK22}
Skopenkov, M. and Krasauskas, R. Surfaces containing two circles through each point, arXiv:1512.09062v3. 

\bibitem{S91}Szőke, R. Complex structures on tangent bundles of Riemannian manifolds. {\em Mathematische Annalen}. \textbf{291} pp. 409-428 (1991)

\bibitem{W90}Wiśniewski, J. On a conjecture of Mukai. {\em Manuscripta Mathematica}. \textbf{68}, 135-141 (1990)

\bibitem{Y91}Yang, C. Any Blaschke manifold of the homotopy type of $\bCP$ has the right volume. {\em Pacific J. Math}. \textbf{151}, 379-394 (1991)

\end{thebibliography}
\end{document}